\documentclass[12pt]{article}
  \usepackage{amsfonts}
\usepackage[centertags]{amsmath}
\usepackage{amssymb}
\usepackage{graphicx}
\newcommand{\ad}{\mathop{\rm ad}\nolimits}

\newcommand{\diag}{\mathop{\rm diag}\nolimits}
\newtheorem{theorem}{Theorem}[section]

\newtheorem{proposition}[theorem]{Proposition}
\newtheorem{definition}[theorem]{Definition}

\newcommand{\Lin}{\mathop{\rm Lin}\nolimits}

\usepackage[verbose,colorlinks=true,naturalnames=true,linkcolor=blue,]{hyperref}

\newcounter{rown}

\begin{document}
\title{{\small\begin{flushright}
In blessed memory of \phantom{aa}\\
Prof. Yu.F. Smirnov \phantom{aaa}\\
(26.05.1935 -- 26.07.2008)\\
\phantom{aaaa}
\end{flushright}}
Extremal projectors for contragredient Lie (super)symmetries (short review)
}

\author{V.N. Tolstoy \\ \\
Lomonosov Moscow State University\\
Skobeltsyn Institute of Nuclear Physics
\\ 
Moscow 119991, Russian Federation}

\date{}
\maketitle \begin{abstract} A brief review of the extremal projectors for contragredient
Lie (super)\-symmetries (finite-dimen\-sional simple Lie algebras, basic classical  Lie
superalgebras, infinite-dimensional affine Kac-Moody algebras and superalgebras, as well
as their quantum $q$-\-analogs) is given. Some bibliographic comments on the applications
of extremal projectors are presented.
\end{abstract}

\setcounter{equation}{0}
\section{Introduction}

Let $G$ be a finite (or compact) group, and $T$ be its unitary representation in a linear
space $V$ -- that is, we have $g\mapsto T(g)$ $(g\in G)$, where $T(g)$ are linear
operators acting in $V$ and satisfying the condition $T(g_1g_2)=T(g_1)T(g_2)$. The
representation $T$ in $V$ is referred to as an irreducible one if $\Lin\{T(G)v\} =V$ for
any nonzero vector $v\in V$. An irreducible representation (IR) is labeled with an
additional superscript $\lambda$, $T^{\lambda}(g)$ (correspondingly, $V^{\lambda}$) or,
in matrix form, $T^{\lambda}(g)=(t^{\lambda}_{ij}(g))$ $(i,j=1,2,\ldots, n)$, where $n$
is the dimension of $V^{\lambda}$.

It is well-known that the elements 
\begin{equation}\label{fg1}
P_{ij}^{\lambda}=\frac{n}{\dim G}\,\sum_{g\in G}T(g)t^{\lambda}_{ij}(g)
\end{equation}
are projection operators for the finite group $G$; that is, they satisfy the following
properties: 
\begin{eqnarray}\label{fg2}
P_{ij}^{\lambda}P_{kl}^{\lambda'}\!\!&=\!\!&\delta_{\lambda\lambda'} \delta_{jk}
P_{il}^{\lambda},\qquad (P_{ij}^{\lambda})^*\,=\,P_{ji}^{\lambda},
\end{eqnarray}
where $"*"$ is Hermitian conjugation.

If $G$ is  a compact group then the formula (\ref{fg1}) is modified as follows: 
\begin{equation}\label{fg3}
P_{ij}^{\lambda}\,=\,\frac{\dim\lambda}{|G|}\int\limits_{g\in G}T(g)
t^{\lambda}_{ij}(g)dg,
\end{equation}
where $|G|$ is the volume of the group $G$ and $dg$ is a $g$-invariant measure on the
group $G$. In the case of $G=SO(3)$ (or $SU(2)$) we have
\begin{equation}\label{fg4}
P_{mm'}^{j}=\frac{2j+1}{8\pi^2}\int_{}T(\alpha,\beta,\gamma)
D^j_{mm'}(\alpha,\beta,\gamma)\sin\beta\,d\alpha\, d\beta\, d\gamma,
\end{equation}
where $\alpha,\beta,\gamma$ are the Euler angles and $D^j_{mm'}(\alpha,\beta,\gamma)$ is
the Wigner $D$-function. The projection operator $P^j:=P^j_{jj}$ is referred to as the
operator of projection onto the highest weight $j$.

Thus, we see that the projection operators in the form (\ref{fg1}) or (\ref{fg3}) require
knowledge of explicit expressions for the operator function $T(g)$ and  the matrix
elements of irreducible representations, $t^{\lambda}_{ij}(g)$, as well as (in the case
of a compact group) the $g$-invariant measure $dg$. In the case of an arbitrary compact
group $G$, the derivation of these expressions involves some problems. There naturally
arises the question of whether it is possible to construct projection operators not in
terms of the elements of a compact group and its representations but in the terms of its
Lie algebra, since the compact group is completely determined by its Lie algebra. The
answer to this question appears to be positive, and the history of its derivation is
rather instructive. We briefly remind this history. It opens with the angular momentum
Lie algebra.

The angular momentum Lie algebra $\mathfrak{so}(3)$ ($\simeq \mathfrak{su}(2)$) is
generated by the three elements (generators) $J_{+}$, $J_{-}$ and $J_{0}$ with the
defining relations: 
\begin{equation}\label{ls1}
[J_0,J_{\pm}]\,=\,\pm J_{\pm},\quad [J_{+},J_{-}]\,=\,2J_0,\quad J_{\pm}^*\,=
\,J_{\mp},\quad J_0^*\,=\,J_0.
\end{equation}
It is obvious that the projector $P^j=P^j_{jj}$, (\ref{fg4}), onto the highest weight $j$
satisfies the relations: 
\begin{equation}\label{ls8}
J_{+}^{}P^j\,=\,P^jJ_{-}^{}\,=\,0,\qquad(P^j)^2\,=\,P^j.
\end{equation}

An associative polynomial algebra of the generators $J_{\pm}$, $J_0$ is called the
universal enveloping algebra of the angular momentum Lie algebra and it is denoted by
$U(\mathfrak{so}(3))$ (or $U(\mathfrak{su}(2))$.

The following proposition holds (no-go theorem): {\it No nontrivial solution of the set
of the equations} \begin{equation} J_{+}^{}P\,=\,PJ_{-}^{}\,=\,0 \label{ls9}
\end{equation}
{\it exists in the space $U(\mathfrak{su}(2))$; that is, the unique solution of this set
of the equations for $P\in U(\mathfrak{su}(2))$ is a trivial one, $P\equiv0$}\footnote{~A
rigorous mathematical formulation of this statement reads as follows: {\it The universal
enveloping algebra $U(\mathfrak{g})$ of a simple Lie algebra $\mathfrak{g}$ has no zero
divisors.}}.

Thus, the theorem states that the projector $P^j$ does not exist in the form of a
polynomial in the generators $J_{\pm}$, $J_0$. This no-go theorem was well known to
mathematicians, but we can assume that it was not known to the majority of physicists.

In 1964 (more than 45 years ago), the Swedish physicist and chemist
P.-O.~L\"ow\-din\,\footnote{~Per-Olov L\"ow\-din was born in 1916 in Uppsala, Sweden, and
died in 2000 (see http://\linebreak[0]www.quantum-chemistry-history.com/Lowdin1.htm).},
who probably did not know the no-go theorem, published a paper in 
{\it Reviews Modern Physics} \cite{L}, where he considered the following operator: 
\begin{equation}\label{ls10}
P^j\,:=\,\prod_{j'\neq j}\frac{{\bf J}^2-j'(j'+1)}{j(j+1)-j'(j'+1)}.
\end{equation}
Here  ${\bf J}^2$ is the Casimir element 
of the angular momentum Lie algebra: 
\begin{equation} {\bf J}^2 \,=\,
\frac{1}{2}\Bigl(J_+^{}J_-^{}
+J_-^{}J_+^{}\Bigr)+J_0^2\,=\,J_-^{}J_+^{}+J_0^{}(J_0^{}+1). \label{ls2}
\end{equation}
The operator (\ref{ls10}) satisfies the equations 
\begin{equation}\label{ls14}
[J_{0}^{},\,P^j]\,=\,0,\quad J_{+}^{}P^j\,=\,P^jJ_{-}^{}\,=\,0,\quad (P^j)^2=P^j
\end{equation}
under the condition that the left- and right-hand sides of these equalities are applied
to vectors characterized by the angular-momentum projection $j$, $J_0\Psi_j=j\Psi_j$.
Therefore, the element (\ref{ls10}) is the operator of projection onto the highest weight
$j$.

After quite involved complicated calculations, L\"owdin reduced the operator (\ref{ls10})
to the form 
\begin{equation}\label{ls16}
P^j\,=\,\sum_{n\ge0}\,\frac{(-1)^n(2j+1)!}{n!(2j+1+n)!}\,J_-^nJ_+^n.
\end{equation}

One year later, in 1965, another physicist, J.~Shapiro from United States of America,
published an article in {\it Journal of Mathematical Physics} \cite{Sh}, where he
proposed the following: "Let us forget the initial expression in the form of infinite
product, (\ref{ls10}), and consider the defining equations (\ref{ls14}), where $P^j$ has
the following ansatz: 
\begin{equation}\label{ls17}
P^j\,=\,\sum_{n\ge0}\,C_n(j)\,J_-^nJ_+^n.\;"
\end{equation}
Substituting this expression into (\ref{ls14}), we immediately obtain the formula
(\ref{ls16}).

It is convenient to get rid of the superscript $j$ in $P^j$ by making the substitution
$j\rightarrow J_0$. As a result, we arrive at 
\begin{eqnarray}\label{ls18}
P&\!\!=&\!\!\sum_{n\ge0}\,\frac{(-1)^n}{n!}\;\varphi_{n}(J_0)\,J_-^nJ_+^n,
\end{eqnarray}
where
\begin{eqnarray}\label{ls19}
\varphi_{n}(J_0)&\!\!=&\!\!\prod_{k=1}^{n}(2J_0+1+k)^{-1}.
\end{eqnarray}
The element $P$ is called the extremal projector. Acting by the extremal projector $P$ on
any weight $\mathfrak{su}(2)$-module $M$ we obtain a space $M^{0}=pM$ of highest weight
vectors for $M$ (if $pM$ has no singularities).

The extremal projector $P$ does not belong to $U(\mathfrak{su}(2))$, but it belongs to
some extension of the universal enveloping algebra. This extension is defined as follows.

Let us consider the formal Taylor series 
\begin{equation}\label{ls22}
\sum_{n,k\ge0}C_{n,k}(J_0)\,J_-^nJ_+^k
\end{equation}
under the condition that $|n-k|\le N$ for some $N\in \mathbb{Z}_+$. The coefficients
$C_{n,k}(J_0)$ are rational functions of the generator $J_0$.

Let $TU(\mathfrak{su}(2))$ be a linear space of such formal series. One can show  that
{\it $TU(\mathfrak{su}(2))$ is an associative algebra with respect to the multiplication
of formal series}. The associative algebra $TU(\mathfrak{su}(2))$ is called the {\it
Taylor extension} of $U(\mathfrak{su}(2))$. It is obvious that $TU(\mathfrak{su}(2))$
contains $U(\mathfrak{su}(2))$.

The extremal projector (\ref{ls18}) belongs to the Taylor extension
$TU(\mathfrak{su}(2))$. Thus, L\"owdin and Shapiro found a solution of the equations
(\ref{ls14}) in the extension of the space $U(\mathfrak{su}(2))$ —- that is, in
$TU(\mathfrak{su}(2))$ -—rather than in the space $U(\mathfrak{su}(2))$ itself.

Later, Shapiro tried to generalized the formula (\ref{ls18}) to the case of
$\mathfrak{su}(3)$ ($\mathfrak{u}(3)$). The Lie algebra $\mathfrak{u}(3)$ is generated by
nine elements $e_{ik}$ $(i,k=1,2,3)$ with the relations: 
\begin{equation}\label{ls26}
[e_{ij},\,e_{kl}]\,=\,\delta_{jk}e_{il}-\delta_{il}e_{kj},\qquad e_{ij}^*\,=\,e_{ij}.
\end{equation}
Shapiro considered an ansatz for $P:=P(\mathfrak{su}(3))$ in the form 
\begin{equation}\label{ls27}
P:=\sum_{\{n_i\},\{m_i\}}C_{\{n_{i}\},\{m_{i}\}}(e_{11},e_{22},e_{33})
\,e_{21}^{n_1}e_{31}^{n_2}e_{32}^{n_3}e_{12}^{m_1} e_{13}^{m_2}e_{23}^{m_3}
\end{equation}
and used the equations 
\begin{eqnarray}\label{ls28}
e_{ij}^{}P&\!\!=&\!\!Pe_{ji}^{}\,=\,0\quad (i<j),\qquad[e_{ii},\,P]\;=\;0\quad (i=1,2,3).
\end{eqnarray}
From the last equations, it follows that 
\begin{equation}\label{ls30}
n_1+n_2-m_1-m_2\,=\,0,\qquad n_2+n_3-m_2-m_3\,=\,0.
\end{equation}
Under the conditions (\ref{ls30}) the expression (\ref{ls27}) belongs to
$TU(\mathfrak{su}(3))$. A set of equations for the coefficients $C_{\{n_{i}\},\{m_{i}\}}
(e_{11},e_{22},e_{33})$ proved to be rather complicated, and Shapiro failed to solve it
and stopped at this stage.

In 1968 R.M.~Asherova and Yu.F.~Smirnov \cite{AS1} made the first important step towards
deriving an explicit formula for the extremal projector for $\mathfrak{u}(3)$. They
proposed applying the Shapiro ansatz (\ref{ls27}) to the extremal projector for the
$\mathfrak{su}(2)$ ($T$-spin) subalgebra generated by the elements $e_{23},\;e_{32}$ and
$e_{22}-e_{33}$; that is \begin{eqnarray}\label{ls31}
P_{23}&\!\!=&\!\!\sum_{n\ge0}\,\frac{(-1)^n}{n!}\;
\varphi_{n}(e_{22}-e_{33})\,e_{32}^ne_{23}^n,
\end{eqnarray}
where \begin{eqnarray}\label{ls32}
\varphi_{n}(e_{22}-e_{33})&\!\!=\!\!&\prod_{k=1}^{n}(e_{22}-e_{33}+1+k)^{-1}.
\end{eqnarray}
Since $e_{23}P_{23}=0$, the operator $P$ takes the form
\begin{equation}\label{ls33}
P=\sum_{n_{i}\ge0}C_{n_{1},n_{2},n_{3}}(e_{11},e_{22},e_{33})
\,e_{21}^{n_1}e_{31}^{n_2}e_{32}^{n_3}e_{12}^{n_1-n_3} e_{13}^{n_2+n_3}\,P_{23}~.
\end{equation}
In this case, the set of equations for the coefficients $C_{\{n_i\}}(\{e_{ii}\})$ becomes
simpler then the corresponding set of equations for the Shapiro ansatz (\ref{ls27}). It
was solved, although the resulting expressions for the coefficients
$C_{\{n_i\}}(\{e_{ii}\})$ are rather complicated.

The next simple idea \cite{T} was to apply, to the expression (\ref{ls33}) from the left,
the extremal projector of the $\mathfrak{su}(2)$ subalgebra generated by
$e_{12},\,e_{21}$ and $e_{11}-e_{22}$. As a result, we obtain the extremal projector
$P(\mathfrak{su}(3))$ in the simple form 
\begin{equation}\label{ls34}
P(\mathfrak{su}(3))=P_{12}\Bigl(\sum_{n\ge0}C_{n}(e_{11}-e_{33})
\,e_{31}^{n}e_{13}^{n}\Bigr)P_{23}.
\end{equation}
The final expression is given by
\begin{equation}\label{ls35}
P(\mathfrak{su}(3))=P_{12}P_{13}P_{23},
\end{equation}
where
\begin{eqnarray}\label{ls36}
P_{ij}&\!\!=&\!\!\sum_{n\ge0}\,\frac{(-1)^n}{n!}\;\varphi_{ij,n}(e_{ii}-e_{jj})\,
e_{ji}^ne_{ij}^n\quad (i<j),
\\[5pt]\label{ls37}
\varphi_{ij,n}(e_{ii}-e_{jj})&\!\!=&\!\! \prod_{k=1}^{n}(e_{ii}-e_{jj}+j-i+k)^{-1}.
\end{eqnarray}
It turned out this formula is a key one. We rewrite this expression for
$P(\mathfrak{su}(3))$ in the terms of the Cartan--Weyl basis with Greek indexes, namely
we replace the root indexes $12$, $23$, $13$ by $\alpha$, $\beta$, $\alpha+\beta$
correspondingly. Moreover we set $h_{\alpha}:=e_{11}-e_{22}$, $h_{\beta}:=e_{22}-e_{33}$,
$h_{\alpha+\beta}:=e_{11}-e_{33}$. In these terms the extremal projector
$P(\mathfrak{su}(3))$ has the form 
\begin{equation}\label{ep3}
P(\mathfrak{su}(3))=P_{\alpha}P_{\alpha+\beta}P_{\beta},
\end{equation}
where
\begin{equation}\label{ep4}
P_{\gamma}^{}=\sum_{n\geq0}\frac{(-1)^{n}} {n!}\,\varphi_{\gamma,n}^{}e_{-\gamma}^{n}\,
e_{\gamma}^{n},
\end{equation}
\begin{equation}\label{ep5}
\varphi_{\gamma,n}^{}=\prod\limits_{k=1}^{n}\Bigl(h_{\gamma}^{}+\rho(\gamma)+
\frac{1}{2}(\gamma,\gamma)k\Bigr)^{-1}.
\end{equation}
Here $\rho$ is a linear function on the positive root system $\Delta_+=
\{\alpha,\beta,\alpha+\beta\}$, such that $\rho(\pi)=\frac{1}{2}(\pi,\pi)$ for all simple
roots $\pi\in\Pi:=\{\alpha,\beta\}$.

Later, the explicit formulas (\ref{ep3})--(\ref{ep5}) were generalized to all
finite-dimensional simple Lie algebras \cite{AS2}--\cite{AST2}, basic classical Lie
superalgebras \cite{T1}, infinite-dimensional affine Kac-Moody algebras and superalgebras
\cite{T2}, as well as all their quantum $q$-analogs \cite{T3,KT1}. At the present time,
the extremal projector method -- that uses explicit formulas of the extremal projectors
-- is a powerful and universal method for solving many problems in the representation
theory, which are widely applied in the theoretical and mathematical physics. For
instance, the method makes it possible to classify irreducible modules, to decompose them
into submodules (for example, to analyze the structure of Verma modules), to describe
reduced (super)algebras (which are associated with the reduction of a (super)algebra to a
subalgebra), to construct bases of representations (for example, the Gelfand-Tsetlin's
type), to develop a detailed theory of Clebsch-Gordan coefficients and other elements of
Wigner-Racah calculus (including compact analytic formulas for these elements and their
symmetry properties), and so on.

The ensuing exposition is organized as follows. In Section 2 we describe the extremal
projectors for finite-dimensional simple Lie algebras. In Section 3 we present the
extremal projectors for finite-dimensional basic classical Lie superalgebras. In Section
4 the extremal projectors for infinite-dimensional affine Kac--Moody algebras and
superalgebras are considered. In Section 5 we describe the extremal projectors for
quantum $q$-analogs of Lie algebras and superalgebras. Finally, in Section 6 some
bibliographic comments on the applications of extremal projectors are presented.

It should be noted that a short review of the extremal projector method was already
published in \cite{T4a}, however there we considered explicit formulas of the extremal
projectors only for simple Lie algebras and we demonstrated some applications of these
projectors for derivation of the general analytic formulas of Clebsch-Gordan coefficients
for the $su(2)$ and $su(3)$ Lie algebras.

\setcounter{equation}{0}
\section{Extremal projectors for simple Lie algebras}

Let $\mathfrak{g}$ be a finite-dimensional simple Lie algebra of rank $r$ and
$\Pi:=\{\alpha_{1},\alpha_{2},\ldots,\alpha_{r}\}$ be its simple root system. Let
$\Delta_+(\mathfrak{g})$ be a system of all positive roots of $\mathfrak{g}$. Any root
$\gamma$ of $\Delta_+(\mathfrak{g})$ has the form: $\gamma=\sum_{i}^{r}l_{i}^{(\gamma)}
\alpha_i$, where all $l_{i}^{(\gamma)}$ are nonnegative integers. A generalization of the
formulas (\ref{ep3})--(\ref{ep5}) to the case of arbitrary finite-dimensional simple Lie
algebra $\mathfrak{g}$ is associated with the concept of normal ordering in the system
$\Delta_+(\mathfrak{g})$.
\begin{definition}\label{df1}
We say that the system $\Delta_+(\mathfrak{g})$ is in normal (convex) ordering if each
composite (not simple) root $\gamma\!=\!\alpha\!+\!\beta$ ($\alpha,\beta,\gamma\in
\Delta_+(\mathfrak{g})$) is written between its components $\alpha$ and $\beta$. It means
that in the normal ordering system $\Delta_+(\mathfrak{g})$ we have either
\begin{eqnarray}\label{CS1}
\mbox{\small$\ldots,\alpha,\ldots,\alpha\!+\!\beta,\ldots,\beta\ldots$},\qquad
\end{eqnarray}
or 
\begin{eqnarray}\label{CS2}
\mbox{\small$\ldots,\beta,\ldots,\alpha\!+\!\beta,\ldots,\alpha,\ldots$}.
\end{eqnarray}
We say also that $\alpha\prec\beta$ if $\alpha$ is located on the left side of $\beta$ in
the normal ordering system $\Delta_+(\mathfrak{g})$, i.e. this corresponds to the case
(\ref{CS1}).
\end{definition}
The normally ordered system $\Delta_+(\mathfrak{g})$ is denoted by the symbol
$\vec{\Delta}_+(\mathfrak{g})$. It is evident that two boundary (end) roots in
$\vec{\Delta}_+(\mathfrak{g})$ are simple.

Let us write down the normal orderings for all simple Lie algebras of rank 2 (see
\cite{AST2}):
\begin{eqnarray}\label{CS3} &A_1\otimes
A_1:\;\mbox{\small$\alpha,\beta\,\leftrightarrow\,\beta,\alpha $},
\\
\label{CS4}&A_2:\;\mbox{\small$\alpha,\alpha\!+\!\beta,\beta\,\leftrightarrow\,
\beta,\alpha\!+\!\beta,\alpha$},
\\
\label{CS5}&B_2:\;\mbox{\small$\alpha,\alpha\!+\!\beta,\alpha\!+\!2\beta,\beta\,
\leftrightarrow\,\beta,\alpha\!+\!2\beta,\alpha\!+\!\beta,\alpha$},
\\
\label{CS6} &G_2:\;\mbox{\small$\alpha,\alpha\!+\!\beta,2\alpha\!+\!3\beta,
\alpha\!+\!2\beta,\alpha\!+\!3\beta,\beta\,\leftrightarrow\,$}
\mbox{\small$\beta,\alpha\!+\!3\beta,\alpha\!+\!2\beta,2\alpha\!+\!3\beta,
\alpha\!+\!\beta,\alpha$},
\end{eqnarray}
where $\alpha-\beta$ is not any root.

{\it We say that a positive root $\gamma\in\Delta_+(\mathfrak{g})$ is generated by
positive roots $\alpha$ and $\beta$ if it is presented in the form of their linear
combination: $\gamma=k\alpha+l\beta$, where $k,l$ are nonzero integers.}

Let $\alpha$ and $\beta$ be any two roots from $\Delta_+(\mathfrak{g})$. We denote by
$\{\alpha,\beta\}$ the subset which contains the roots $\alpha$, $\beta$ and all positive
roots from $\Delta_+(\mathfrak{g})$, generated by the roots $\alpha$ and $\beta$.
\begin{definition}
Let in the normal ordering system $\vec{\Delta}_+(\mathfrak{g})$ between two roots from
the subset $\{\alpha,\beta\}$ there are not another roots except the roots generated by
themselves $\alpha$ and $\beta$. Then inverting this subset we again obtain the normal
ordering system $\vec{\Delta}'_+(\mathfrak{g})$. Such transformation is called the
elementary inversion.
\end{definition}
Because any two positive roots generate a root system of a Lie algebra of rank 2,
therefore all elementary inversions coincide with the elementary inversions of normal
orders for the root systems $\Delta_+(\mathfrak{g})$ of the Lie algebras $\mathfrak{g}$
of rank 2, (\ref{CS3})--(\ref{CS6}).

The combinatorial structure of the root system $\Delta_+(\mathfrak{g})$ is described the
following theorem. 
\begin{theorem}\label{th1}
(i) Normal ordering in the system $\Delta_+(\mathfrak{g})$ exists for any mutual location
of the simple roots $\alpha_i,\; i=1,2,\ldots,r$. (ii) Any two normal orderings
$\vec{\Delta}_+(\mathfrak{g})$ and $\vec{\Delta}'_+ (\mathfrak{g})$ can be obtained one
from another by compositions of elementary inversions for the root systems of the Lie
algebras of rank 2.
\end{theorem}
A detailed proof of the theorem is presented in \cite{AST2, T4}.

Let $e_{\pm\gamma}^{}$, $h_{\gamma}^{}$ be Cartan-Weyl root vectors normalized by the
condition 
\begin{equation}\label{ep1}
[e_{\gamma}^{},\,e_{-\gamma}^{}]=h_{\gamma}^{}.
\end{equation}
We construct a formal Taylor series on the following monomials
\begin{eqnarray}\label{TE2} e_{-\beta}^{n_{\beta}}\cdots
e_{-\gamma}^{n_{\gamma}}e_{-\alpha}^{n_{\alpha}}\: e_{\alpha}^{m_{\alpha}}
e_{\gamma}^{m_{\gamma}}\cdots e_{\beta}^{m_{\beta}}
\end{eqnarray}
with coefficients which are rational functions of the Cartan elements $h_{\alpha_i}$
($i=1,2,\ldots,r$), and nonnegative integers $n_{\beta},\ldots,n_\gamma,n_{\alpha},
m_{\alpha},m_\gamma,\ldots,m_{\beta}$ are subjected to the constraints (for some $N\in
\mathbb{Z}_+$) 
\begin{eqnarray}\label{TE3}
\Bigr|\sum_{\gamma\in\Delta_{+}}(n_{\gamma}-m_{\gamma})l_{i}^{(\gamma)}\Bigl|\leq N
, \quad i=1,2,\cdots,r,
\end{eqnarray}
for all monomial of the given series. Here $l_{i}^{(\gamma)}$ are coefficients in a
decomposition of the root $\gamma$ with respect to the system of simple roots $\Pi$,
$\gamma=\sum_{i}^{r}l_{i}^{(\gamma)}\alpha_i$, $\alpha_i\in\Pi$. Let $TU(\mathfrak{g})$
be a linear space of all such formal series. We have the following simple proposition.
\begin{proposition}\label{th2} The linear space $TU(\mathfrak{g})$ is an associative
algebra with respect to a multiplication of formal series.
\end{proposition}
The algebra $TU(\mathfrak{g})$ is called the Taylor extension of $U(\mathfrak{g})$.
\begin{theorem}\label{th3} The equations 
\begin{equation}\label{ep2}
e_{\gamma}^{}P(\mathfrak{g})=P(\mathfrak{g})e_{-\gamma}^{}=0\quad(\forall\;
\gamma\in\Delta_+(\mathfrak{g}))~, \qquad P^{2}(\mathfrak{g})=P(\mathfrak{g})
\end{equation}
have a unique nonzero solution in the space $TU(\mathfrak{g})$, $P(\mathfrak{g})\in
TU(\mathfrak{g})$, and this solution has the form 
\begin{equation}\label{ep23}
P(\mathfrak{g})\,=\,\prod_{\gamma\in\vec{\Delta}_+(\mathfrak{g})}P_{\gamma}^{}
\end{equation}
for any normal ordering system $\vec{\Delta}_+(\mathfrak{g})$, where the elements
$P_{\gamma}$ are defined by the formulae
\begin{equation}\label{ep24}
P_{\gamma}^{}\,=\,\sum_{n\geq0}\frac{(-1)^{n}}{n!}\,\varphi_{\gamma,n}^{}
e_{-\gamma}^{n}\,e_{\gamma}^{n},
\end{equation}
\begin{equation}\label{ep25}
\varphi_{\gamma,n}^{}\,=\,\prod\limits_{k=1}^{n}\Bigl(h_{\gamma}^{}+\rho(\gamma)+
\frac{1}{2}(\gamma,\gamma)k\Bigr)^{-1},
\end{equation}
and $\rho$ is the linear function on the positive root system $\Delta_+(\mathfrak{g})$,
such that $\rho(\alpha_i)=\frac{1}{2}(\alpha_i,\alpha_i)$ for all simple roots
$\alpha_i\in\Pi$\footnote{In the case of the finite-dimensional simple Lie algebras, the
function $\rho(\cdot)$ can be presented in the form of the scalar product
$\rho(\gamma)=(\rho,\gamma)$, where $\rho$ on the right-hand side is a half-sum of all
positive roots.}. Thus, the extremal projector $P(\mathfrak{g})$ does not depend on
choice of normal order in $\Delta_+(\mathfrak{g})$.
\end{theorem}
The extremal projectors of the Lie algebras of rank 2 and the combinatorial theorem play
key roles for the proof of this theorem for an arbitrary simple Lie algebra
$\mathfrak{g}$ of rank $r>2$ \cite{AST2, T4}. The extremal projectors of the Lie algebras
of rank 2 are given by (see \cite{AST2})
\begin{eqnarray}\label{CS7}
P_{\alpha}P_{\beta}\!\!&=\!\!&P_{\beta}P_{\alpha}\qquad\; (A_1\otimes A_1),
\\[3pt]
\label{CS8}P_{\alpha}P_{\alpha\!+\!\beta}P_{\beta}\!\!&=\!\!&
P_{\beta}P_{\alpha\!+\!\beta}P_{\alpha}\qquad\qquad(A_2),
\\[3pt]
\label{CS9}P_{\alpha}P_{\alpha\!+\!\beta}P_{\alpha\!+\!2\beta,\beta} \!\!&=\!\!&
P_{\beta}P_{\alpha\!+\!2\beta}P_{\alpha\!+\!\beta}P_{\alpha}\qquad\quad(B_2),
\\[3pt]
\label{CS10}P_{\alpha}P_{\alpha\!+\!\beta}P_{2\alpha\!+\!3\beta} P_{\alpha\!+\!2\beta}
P_{\alpha\!+\!3\beta}P_{\beta}\!\!&=\!\!& P_{\beta}P_{\alpha\!+\!3\beta}P_{\alpha\!+
\!2\beta}P_{2\alpha\!+\!3\beta} P_{\alpha\!+\!\beta}P_{\alpha}\quad(G_2).
\end{eqnarray}
We can show that these equations are valid not only for $\rho(\alpha)=\frac{1}{2}
(\alpha,\alpha)$ and $\rho(\beta)=\frac{1}{2}(\beta,\beta)$ but as well as for
$\rho(\alpha)=x_\alpha$ and $\rho(\beta)=x_\beta$ where $x_\alpha$ and $x_\beta$ are
arbitrary complex numbers. Taking into account the second part {\it(ii)} of the theorem
(\ref{th1}) it now follows the equalities (\ref{ep2}) for (\ref{ep23})--(\ref{ep25})
immediately. A uniqueness of this solution in the space $TU(\mathfrak{g})$ is proved more
complicated \cite{Zh}.

\setcounter{equation}{0}
\section{Extremal projectors for Lie superalgebras}

In this section we consider the extremal projectors for finite-dimensional basic
classical Lie superalgebras. Let $\mathfrak{g}$ be a finite-dimensional basic classical
Lie superalgebra \cite{Kac1, Kac2} of rank $r$ and $\Pi:=\{\alpha_{1},
\alpha_{2},\ldots,\alpha_{r}\}$ be its simple root system. Let $\Delta_+(\mathfrak{g})$
be a system of all positive roots of $\mathfrak{g}$. Any root $\gamma$ of
$\Delta_+(\mathfrak{g})$ has the form: $\gamma=\sum_{i}^{r} l_i^{(\gamma)}\alpha_i$,
where all $l_i^{(\gamma)}$ are nonnegative integers. There are two type of the roots:
even and odd, $\Delta_+(\mathfrak{g})= \Delta_{+}^{(0)}
(\mathfrak{g})\oplus\Delta_{+}^{(1)}(\mathfrak{g})$, and they are classified (colored) as
follows: 
\begin{itemize}
\item
Any even root $\gamma\in\Delta_+(\mathfrak{g}^{(0)})$ is white. In this case
$2\gamma\notin \Delta_{+}^{(0)}(\mathfrak{g})$ and $(\gamma,\gamma)\ne0$.
\item
A odd root $\gamma\in\Delta_{+}^{(1)}(\mathfrak{g})$ is called grey if $2\gamma\notin
\Delta_+(\mathfrak{g}^{(0)})$. In this case $(\gamma,\gamma)=0$.
\item
A odd root $\gamma\in\Delta_{+}^{(1)}(\mathfrak{g})$ is called dark if $2\gamma\in
\Delta_{+}^{(0)}(\mathfrak{g})$. In this case $(\gamma,\gamma)\ne0$.
\end{itemize}
The total system of all roots, $\Delta(\mathfrak{g})$, has the form:
$\Delta(\mathfrak{g})=\Delta(\mathfrak{g})_+\bigcup(-\Delta_+(\mathfrak{g}))$ and the
parity of the negative root -$\gamma$ coincides with the parity of the positive root
$\gamma$.

Now we define of a reduced system of the positive root system $\Delta_+(\mathfrak{g})$.
\begin{definition}\label{df3} 
The reduced system is called a set $\underline{\Delta}_+(\mathfrak{g})$ obtained from the
positive system $\Delta_+(\mathfrak{g})$ by removing of all doubled roots $2\gamma$ where
$\gamma$ is a dark odd root, that is $\underline{\Delta}_+(\mathfrak{g})= \Delta_+
(\mathfrak{g})\backslash\{2\gamma\in\Delta_+(\mathfrak{g})|\,\gamma\;{\rm is\; odd}\}.
$\footnote{\rm In the case of any simple Lie algebra $\mathfrak{g}$, the reduced system
coincides with the total positive root system, $\underline{\Delta}_+(\mathfrak{g})
=\Delta_+(\mathfrak{g})$.}
\end{definition}
Normal ordering in the  system $\underline{\Delta}_+(\mathfrak{g})$ is fixed by the
definition (\ref{df1}). The combinatorial theorem (\ref{th1}) is also true for the
reduced system $\underline{\Delta}_+(\mathfrak{g})$. In this case there are the following
elementary inversions of the reduced systems of the Lie algebras and superalgebras of
rank 2 (see \cite{KT2}):
\begin{eqnarray}\label{CS15}
&\mbox{\small$\alpha,\beta\,\leftrightarrow\,\beta,\alpha $},
\\
\label{CS16}
&\mbox{\small$\alpha,\alpha\!+\!\beta,\beta\,\leftrightarrow\,
\beta,\alpha\!+\!\beta,\alpha$},
\\
\label{CS17} &\mbox{\small$\alpha,\alpha\!+\!\beta,\alpha\!+\!2\beta,\beta\,
\leftrightarrow\,\beta,\alpha\!+\!2\beta,\alpha\!+\!\beta,\alpha$},
\\
\label{CS18}&\mbox{\small$\alpha,\alpha\!+\!\beta,2\alpha\!+\!3\beta,
\alpha\!+\!2\beta,\alpha\!+\!3\beta,\beta\,\leftrightarrow\,$}
\mbox{\small$\beta,\alpha\!+\!3\beta,\alpha\!+\!2\beta,2\alpha\!+\!3\beta,
\alpha\!+\!\beta,\alpha$},
\end{eqnarray}
where $\alpha-\beta$ is not any root.\footnote{Either of the two roots $\alpha$ and
$\beta$ in (\ref{CS15}) can be white, gray or dark. Either of the two roots $\alpha$,
$\beta$ in (\ref{CS16}) and also the root $\alpha$ in (\ref{CS17}) is white or gray. The
root $\beta$ in (\ref{CS17}) is white or dark. The roots $\alpha$ and $\beta$ in
(\ref{CS18}) are white.}

The Cartan-Weyl root vectors $e_{\pm\gamma}^{}$, $h_{\gamma}^{}$
($\gamma\in\underline{\Delta}_+(\mathfrak{g})$) are normalized by the condition
\begin{equation} 
\label{ep11} [e_{\gamma}^{},e_{-\gamma}^{}]=h_{\gamma}^{}.
\end{equation}
The Taylor extension $TU(\mathfrak{g})$ of $U_(\mathfrak{g})$ is generated by the formal
series on the monomials (\ref{TE2}) where the roots $\alpha$, $\gamma$, $\ldots$, $\beta$
belong to $\underline{\Delta}_+(\mathfrak{g})$.

The extremal projector, namely the element $P(\mathfrak{g})$ satisfying the conditions
(\ref{ep2}), is given by the formula \cite{T1} 
\begin{equation}\label{ep27}
P(\mathfrak{g})\,=\,\prod_{\gamma\in\vec{\underline{\Delta}}_+(\mathfrak{g})}
P_{\gamma}^{},
\end{equation}
where the factors $P_{\gamma}$ are defined as follows. If the root $\gamma$ is white,
then
\begin{equation}\label{ep28}
P_{\gamma}^{}\,=\,\sum_{n\geq0}\frac{(-1)^{n}}{n!}\,\varphi_{\gamma,n}^{}
e_{-\gamma}^{n}\,e_{\gamma}^{n},
\end{equation}
where
\begin{equation}\label{ep29}
\varphi_{\gamma,n}^{}\,=\,\prod\limits_{k=1}^{n}\Bigl(h_{\gamma}^{} +\rho(\gamma)+
\frac{1}{2}(\gamma,\gamma)k\Bigr)^{-1}.
\end{equation}
If the root $\gamma$ is gray, then
\begin{equation}\label{ep30}
P_{\gamma}^{}\,=\,1-\frac{1}{h_{\gamma}^{}+\rho(\gamma)}\,e_{-\gamma}^{}e_{\gamma}^{}.
\end{equation}
If the root $\gamma$ is dark, then the factors $P_{\gamma}$ is given by (cf.
\cite{BerTol} for the $\mathfrak{osp}(1|2)$ case)
\begin{equation}\label{ep28'}
P_{\gamma}^{}\,=\,\sum_{n\geq0}\,\frac{1}{n!}\,\Bigl(\varphi_{\gamma,2n}^{(0)}
e_{-\gamma}^{2n}e_{\gamma}^{2n}+\varphi_{\gamma,2n+1}^{(1)}
e_{-\gamma}^{2n+1}e_{\gamma}^{2n+1}\Bigr),
\end{equation}
where
\begin{equation}\label{ep29'}
\varphi_{\gamma,2n}^{(0)}\,=\,2^{-n}\prod\limits_{k=1}^{n}\Bigl(h_{\gamma}^{}
+\rho(\gamma)+ \frac{1}{2}(\gamma,\gamma)(k-1)\Bigr)^{-1},
\end{equation}
\begin{equation}\label{ep29''}
\varphi_{\gamma,2n+1}^{(1)}\,=\,-2^{-n}\prod\limits_{k=1}^{n+1}\Bigl(h_{\gamma}^{}
+\rho(\gamma)+\frac{1}{2}(\gamma,\gamma)k\Bigr)^{-1}.
\end{equation}
The linear function $\rho(\gamma)$, $\rho(\alpha_i)= \frac{1}{2} (\alpha_i,\alpha_i)$ for
all simple roots $\alpha_i\in\Pi$, can be presented in the form of the scalar product
$\rho(\gamma)=(\rho,\gamma)$, where $\rho$ on the right-hand side is a difference between
a half-sum of all even positive roots and a half-sum of all odd positive roots.

\setcounter{equation}{0}
\section{Extremal projectors for affine Kac--Moody\\
(super)algebras}

Let $\mathfrak{g}$ be a affine Kac--Moody (super)algebra \cite{Kac3,vLeur} of rank $r$
and $\Pi\!=\!\{\alpha_{1},\alpha_{2},\ldots, \alpha_{r}\}$ be its simple root system. Let
$\Delta_+(\mathfrak{g})$ be a system of all positive roots of $\mathfrak{g}$. Any root
$\gamma$ of $\Delta_+(\mathfrak{g})$ has the form: $\gamma=\sum_{i}^{r}
l_{i}^{(\gamma)}\alpha_i$, where all $l_{i}^{(\gamma)}$ are nonnegative integers. In this
case the system $\Delta_+(\mathfrak{g})=\Delta_{+}^{(0)}(\mathfrak{g})\oplus
\Delta_{+}^{(1)}(\mathfrak{g})$ has infinite number of the roots which are categorized
into real and imaginary. Every imaginary root $m\delta$ ($m\in\mathbb{N}$) satisfies the
condition $(m\delta,\gamma)=0$ for all $\gamma\in\Delta(\mathfrak{g}) \bigl(=
\Delta_+(\mathfrak{g})\bigcup(-\Delta_+(\mathfrak{g})\bigr)$. For the real roots this
condition is not valid. We again fix the reduced system $\underline{\Delta}_+
(\mathfrak{g})$, which is obtained from the total system $\Delta_+(\mathfrak{g})$ by
removing of all doubled roots $2\gamma$ where $\gamma$ is a dark odd root, and then we
define the normal ordering in $\underline{\Delta}_+(\mathfrak{g})$ as follows. 
\begin{definition} We say that the system $\underline{\Delta}_+(\mathfrak{g})$
is in normal ordering if its roots are written as follows: (i) all imaginary roots are
immediately adjacent to one another, (ii) each composite (not simple) root
$\gamma=\alpha+\beta$ ($\alpha,\beta, \gamma\in \underline{\Delta}_{+}(\mathfrak{g})$),
where $\alpha$ and $\beta$ are not proportional roots ($\alpha\ne\lambda\beta$), is
written between its components $\alpha$ and $\beta$.
\end{definition}
The combinatorial theorem (\ref{th1}) is also true for the reduced system
$\underline{\Delta}_+(\mathfrak{g})$. In this case there are the following elementary
inversions of the reduced systems of algebras and superalgebras of rank 2 (see \cite{T4,
T5}:
\begin{eqnarray}\label{CS11'}
&\mbox{\small$\alpha,\beta\,\leftrightarrow\,\beta,\alpha $},
\\[5pt]
\label{CS12'}&\mbox{\small$\alpha,\alpha+\beta,\beta\,\leftrightarrow\,
\beta,\alpha+\beta,\alpha$},
\\[5pt]
\label{CS13'}&\mbox{\small$\alpha,\alpha+\beta,\alpha+2\beta,\beta\,
\leftrightarrow\,\beta,\alpha+2\beta,\alpha+\beta,\alpha$}.
\\[5pt]
\label{CS14'}&\mbox{\small$\alpha,\alpha\!+\!\beta,2\alpha\!+\!3\beta,
\alpha\!+\!2\beta,\alpha\!+\!3\beta,\beta\,\leftrightarrow\,$}
\mbox{\small$\beta,\alpha\!+\!3\beta,\alpha\!+\!2\beta,2\alpha\!+\!3\beta,
\alpha\!+\!\beta,\alpha$}.
\end{eqnarray}
\begin{equation}\label{CS15'} \begin{array}{l}
\mbox{\small$\alpha,\delta+\alpha,2\delta+\alpha,\ldots,\infty\delta+
\alpha,\delta,2\delta,3\delta,\ldots,\infty\delta,
\infty\delta-\alpha,\ldots,2\delta-\alpha,\delta-\alpha\,\leftrightarrow$}
\\[3pt]
\mbox{\small$\leftrightarrow\,\delta-\alpha,2\delta-\alpha,\ldots,\infty\delta-
\alpha,\delta,2\delta,3\delta,\ldots,\infty\delta,\infty\delta+\alpha,\ldots,2\delta+
\alpha,\delta+\alpha,\alpha$},
\end{array}
\end{equation}
\begin{equation}\label{CS16'} \begin{array}{l} \!\!\!
\mbox{\small$\alpha,\delta\!+\!2\alpha,\delta\!+\!\alpha,3\delta\!+\!2\alpha,
2\delta\!+\!\alpha,\ldots,\infty\delta\!+\!\alpha,(2\infty\!+\!1)\delta\!+\!2\alpha,
(\infty\!+\!1)\delta\!+\!\alpha,\delta,2\delta,\ldots,$}
\\[3pt]
\!\!\!\mbox{\small$\infty\delta,(\infty\!+\!1)\delta\!-\!\alpha,
(2\infty\!+\!1)\delta\!-\!2\alpha,\infty\delta-\alpha,\ldots,
2\delta\!-\!\alpha,3\delta\!-
\!2\alpha,\delta\!-\!\alpha,\delta\!-\!2\alpha,\,\leftrightarrow$}
\\[5pt]
\!\!\!
\mbox{\small$\leftrightarrow\delta\!-\!2\alpha,\delta\!-\!\alpha,3\delta\!-\!2\alpha,
2\delta\!-\!\alpha,\ldots, \infty\delta\!-\!\alpha,(2\infty\!+\!1)\delta\!-\!2\alpha,
(\infty\!+\!1)\delta\!-\!\alpha, \delta,2\delta,\ldots,$}
\\[3pt]
\!\!\! \mbox{\small$\infty\delta,(\infty\!+\!1)\delta\!+\!\alpha,
(2\infty\!+\!1)\delta\!+\!2\alpha,\infty\delta\!+\!\alpha,\ldots,
2\delta\!+\!\alpha,3\delta\!+\!2\alpha,\delta\!+\!\alpha,\delta\!+\!2\alpha,\alpha$},
\end{array}
\end{equation}
where $\alpha-\beta$ is not any root.

The Cartan-Weyl root vectors $e_{\pm\gamma}^{}$, $h_{\gamma}^{}$
($\gamma\in\underline{\Delta}_+(\mathfrak{g}))$ are normalized by the
condition\footnote{If a root space $\mathfrak{g}_{m\delta}$ ($\mathfrak{g}_{-m\delta}$)
of the imaginary root $m\delta$ (-$m\delta$), $m\in\mathbb{N}$, is not multiplicity free,
$\dim\mathfrak{g}_{m\delta}=\dim\mathfrak{g}_{-m\delta}>1$, then we choose basis root
vectors $\{e_{-m\delta,i}\}$ and $\{e_{m\delta,i}\}$ that are dual with respect to a
standard bilinear form on $\mathfrak{g}$, $(e_{-m\delta,i},e_{m\delta,j})=\delta_{ij}$.
In this case $[e_{m\delta,i},e_{-m\delta,j}]=\delta_{ij}h_{m\delta}$, and moreover
$[e_{\pm m\delta,i},e_{\pm m\delta,j}]=0$.}
\begin{equation}\label{aff6}
[e_{\gamma}^{},e_{-\gamma}^{}]=h_{\gamma}^{}.
\end{equation}
We construct a formal Taylor series on the following monomials
\begin{eqnarray}\label{aff7} e_{-\beta}^{n_{\beta}}\cdots
e_{-\gamma}^{n_{\gamma}}e_{-\alpha}^{n_{\alpha}}\, e_{\alpha'}^{m_{\alpha'}}
e_{\gamma'}^{m_{\gamma'}}\cdots e_{\beta'}^{m_{\beta'}}\quad{\rm (finite\;\; product)}
\end{eqnarray}
with coefficients which are rational functions of the Cartan elements $h_{\alpha_i}$
($i=1,2,\ldots,r$), and nonnegative integers $n_{\beta},\ldots,n_\gamma,n_{\alpha},
m_{\alpha'},m_{\gamma'},\ldots,m_{\beta'}$ are subjected to the constraints 
\begin{eqnarray}\label{aff8}
\Bigr|\sum_{\gamma\in\underline{\Delta}_{+}}n_{\gamma}l_{i}^{(\gamma)}-
\sum_{\gamma'\in\underline{\Delta}_{+}}m_{\gamma'} l_{i}^{(\gamma')}\Bigl|\leq N,\quad
i=1,2,\cdots,r,
\end{eqnarray}
for all monomial of the given series. Here $l_{i}^{(\gamma)}$ are coefficients in a
decomposition of the root $\gamma$ with respect to the system of simple roots $\Pi$. Let
$TU_(\mathfrak{g})$ be a linear space of all such formal series. {\it The linear space
$TU_(\mathfrak{g})$ is an associative algebra with respect to a multiplication of formal
series}. The algebra $TU(\mathfrak{g})$ is called the Taylor extension of
$U(\mathfrak{g})$.

The extremal projector $P(\mathfrak{g})$ of the affine Kac--Moody (super)algebra is given
by the formula \cite{T2}:
\begin{equation}\label{ep27'}
P(\mathfrak{g})\,=\,\prod_{\gamma\in\vec{\underline{\Delta}}_+(\mathfrak{g})}
P_{\gamma}^{},
\end{equation}
where the factors $P_{\gamma}$ are defined by the formulas (\ref{ep28})--(\ref{ep29''})
for all real roots. If the root $\gamma$ is even imaginary, $\gamma=m\delta$ ($m\in
\mathbb{N}$), then\footnote{In the case of the superalgebra $\mathfrak{g}=A^{(4)}(2k,2l)$
the imaginary roots can be white (even) as well as dark (odd). If the root $m\delta$ is a
dark imaginary root, then the formula (\ref{ep29'}) is modified slightly by analogy with
(\ref{ep28'}).}
\begin{equation}\label{ep28"}
P_{m\delta}^{}\,=\,\prod_{i=0}^{\dim\mathfrak{g}_{m\delta}}P_{m\delta}^{(i)}
\end{equation}
\begin{equation}\label{ep29'}
P_{m\delta}^{(i)}\,=\,\sum_{n\geq0}\,\frac{(-1)^{n}}{n!(h_{m\delta}+
\rho(m\delta)\bigr)^{n}}\,e_{-m\delta,i}^{n}e_{m\delta,i}^{n}.
\end{equation}
The linear function $\rho(\gamma)$ is determined by $\rho(\alpha_i)=\frac{1}{2}
(\alpha_i,\alpha_i)$ for all simple roots $\alpha_i\in\Pi$.

\setcounter{equation}{0}
\section{Extremal projectors for quantum Lie (super)\-algebras}

Let $\mathfrak{g}(A,\tau)$ be a contragredient Lie (super)algebra of finite
growth\footnote{These (super)algebras include all finite-dimensional simple Lie algebras
and basic classical superalgebras, infinite-dimensional affine Kac-Moody algebras
\cite{Kac3} and superalgebras \cite{vLeur}.} with a symmetrizable Cartan matrix $A$ (i.e.
$A=DA^{sym}$, where $A^{sym}\!=\!(a_{\;ij}^{sym}\!)_{i,j\in I}$ is a symmetrical matrix,
and $D$ is an invertible diagonal matrix, $D=\diag(d_1,d_2,\ldots,d_r)$), $\tau\subset
I$, $I:=\{1,2,\ldots,r\}$, and let $\Pi:=\{\alpha_{1},\ldots,\alpha_{r}\}$ be a system of
simple roots for $g(A,\tau)$. 
\begin{definition}\cite{T3,KT2} The quantum (super)algebra
$U_{q}(\mathfrak{g})$ (where $\mathfrak{g}:= \mathfrak{g}(A,\tau)$) is an associative
(super)algebra over $\mathbb{C}[q,q^{-1}]$ with Chevalley generators
$e_{\pm\alpha_{i}}^{}$, $k_{\alpha_{i}}^{\pm 1}\!:=\!q^{\pm h_{\alpha_{i}}}$, $(i\in
I\!:=\!\{1,2,\ldots,r\})$, and the defining relations: 
\begin{eqnarray}\label{QA1}
k_{\alpha_{i}}^{}k_{\alpha_{i}}^{-1}\!\!&=\!\!&k_{\alpha_{i}}^{-1}k_{\alpha_{i}}^{}\;=
\;1,\qquad\quad k_{\alpha_{i}}k_{\alpha_{j}}\;=\;k_{\alpha_{j}}k_{\alpha_{i}},
\\[7pt]\label{QA2}
k_{\alpha_{i}}^{}e_{\pm\alpha_{j}}^{}k_{\alpha_{i}}^{-1}\!\!&=\!\!&
q^{\pm(\alpha_{i},\alpha_{j})}e_{\pm\alpha_{j}}^{},\qquad[e_{\alpha_{i}},
e_{-\alpha_{j}}] =\delta_{ij}\frac{k_{\alpha_{i}}^{}-k_{\alpha_{i}}^{-1}}{q-q^{-1}},
\\[7pt]\label{QA3}
(\ad_{q}e_{\pm\alpha_{i}})^{n_{ij}}e_{\pm\alpha_j}\!\!&=\!\!&0\qquad{\rm for}\;\;i\neq j,
\end{eqnarray}
where the positive integers $n_{ij}$ are given as follows: $n_{ij}=1$ if
$a_{ii}^{sym}=a_{ij}^{sym}=0$, $n_{ij}=2$ if $a_{ii}^{sym}=0,\;a_{ij}^{sym}\ne0$, and
$n_{ij}=-2a_{ij}^{sym}/a_{ii}^{sym}+1$ if $a_{ii}^{sym}\neq0$. Moreover, if any three
simple roots $\alpha_{i}$, $\alpha_{j}$, $\alpha_{k}^{}$ satisfy the conditions
$(\alpha_i,\alpha_i)=(\alpha_j,\alpha_k)=0$ and $(\alpha_i,\alpha_j)=-(\alpha_i,\alpha_k)
\ne0$, then there are the additional triple relations of the form\footnote{What we really
consider here is a special case when the system  $\Pi:=\{\alpha_{1},\ldots,\alpha_{r}\}$
has a minimal number of odd roots. In the case if $\Pi$ does not satisfies this condition
then another additional Serre's relations can exist (see \cite{Yamane1,Yamane2}).}:
\begin{equation}\label{QA4}
[[e_{\pm\alpha_{i}},e_{\pm\alpha_{j}}]_q,[e_{\pm\alpha_{i}},e_{\pm\alpha_{k}}]_q]_q=0.
\end{equation}
\end{definition}
Here in (\ref{QA1})--(\ref{QA4}) the brackets $[\cdot,\cdot]$ is the usual
supercommutator and $[\cdot,\cdot]_{q}$ and $\ad_{q}$ denote the $q$-deformed
supercommutator ($q$-supercommutator) in $U_{q}(\mathfrak{g})$ \cite{T3}: 
\begin{eqnarray}\label{QA5}
(\ad_{q}e_{\alpha})e_{\beta}\equiv[e_{\alpha},e_{\beta}]_q\!\!&=\!\!&
e_{\alpha}e_{\beta}-(-1)^{\deg(e_{\alpha})\deg(e_{\beta})}q^{(\alpha,\beta)}
e_{\beta}e_{\alpha},
\end{eqnarray}
where $(\alpha,\beta)$ is a scalar product of the roots $\alpha$ and $\beta$, and the
parity function $\deg(\cdot)$ is given by
\begin{eqnarray}\label{QA6}
\deg(k_{\alpha_{i}})=0\;\,(i\in I),\quad \deg(e_{\pm\alpha_i})=
0\;\,(i\not\in\tau),\quad\deg(e_{\pm\alpha_i})=1\;\,(i\in\tau).
\end{eqnarray}
Below we shall use the following short notation:
\begin{eqnarray}\label{QA7}
\vartheta(\gamma)\!\!&:=\!\!&\deg(e_{\gamma}).
\end{eqnarray}
The definition of a quantum algebra also includes operations of  a co-multiplication
$\Delta_{q}$, an antipode $S_{q}$, and a co-unit $\epsilon_{q}$. Explicit formulas of
these operations will not be used here and they are not given.

The $q$-analog of the Cartan-Weyl basis for $U_{q}(\mathfrak{g})$ is constructed by using
the following inductive algorithm \cite{T3,KT1,KT2}, \cite{KT3} -- \cite{KT5}.

{\it We fix some normal ordering $\vec{\underline{\Delta}}_+(\mathfrak{g})$ and put by
induction
\begin{eqnarray}\label{QA7}
e_{\gamma}:=[e_{\alpha},e_{\beta}]_{q},\qquad\quad
e_{-\gamma}:=[e_{-\beta},e_{-\alpha}]_{{q}^{-1}}
\end{eqnarray}
if $\gamma\!=\!\alpha\!+\!\beta$, $\alpha\prec\gamma\prec\beta$
($\alpha,\beta,\gamma\in\vec{\underline{\Delta}}_+(\mathfrak{g})$), and the segment
$[\alpha;\beta]\subseteq\vec{\underline{\Delta}}_+(\mathfrak{g})$ is minimal one
including the root $\gamma$, i.e. the segment has not another roots $\alpha'$ and
$\beta'$ such that $\alpha'\!+\!\beta'\!=\!\gamma$. Moreover we put 
\begin{eqnarray}\label{QA7}
k_{\gamma}:=\prod_{i=1}^{r}k_{\alpha_{i}}^{l_{i}^{(\gamma)}},
\end{eqnarray}
if $\gamma=\sum_{i=1}^{r}l_{i}^{(\gamma)}\alpha_{i}$
($\gamma\in\vec{\underline{\Delta}}_+(\mathfrak{g}),\;\alpha_i\in\Pi$)}.

\noindent By this procedure one can construct the total quantum Cartan-Weyl basis for all
quantized finite-dimensional simple contragredient Lie (super)algebras. In the case of
the quantized infinite-dimensional affine Kac-Moody (super)algebras we have to apply one
more additional condition. Namely, first we construct all root vectors $e_\gamma$
($\gamma\in\underline{\Delta}_+(\mathfrak{g})$) by means of the given procedure, and then
we overdeterminate the generators $e_{\pm m\delta}$ of the imaginary roots
$m\delta\in\underline{\Delta}_+(\mathfrak{g})$ ($m\in\mathbb{N}$) in a way that the new
generators $e_{m\delta}'$ are mutually commutative if they are not conjugate generators.
Because of the fact that we do not have a sufficient place here to describe the
overdetermination of imaginary root generators in details, we are restricted to a
consideration of finite-dimensional case, i.e. when $\mathfrak{g}$ is a
finite-dimensional simple contragredient Lie (super)algebra.

The quantum Cartan-Weyl basis is characterized by the following properties
\cite{T3,KT1,KT2}, \cite{KT3} -- \cite{KT7}.
\begin{proposition}
The root vectors $\{e_{\pm\gamma}\}$ ($\gamma\,\in\underline{\Delta}_+(\mathfrak{g})$)
satisfy the following relations:
\begin{eqnarray}\label{QA8}
k_{\alpha}^{\pm 1}e_{\gamma}\!\!&=\!\!&q^{\pm(\alpha,\gamma)}e_{\gamma}
k_{\alpha}^{\pm1},
\\[7pt]\label{QA9}
[e_{\gamma}^{},e_{-\gamma}^{}]\!\!&=\!\!&a(\gamma)\frac{k_{\gamma}-k_{\gamma}^{-1}}
{q-q^{-1}},
\\[7pt]\label{QA10}
[e_{\alpha}^{},e_{\beta}^{}]_{q}\!\!&=\!\!&\sum_{\alpha\prec\gamma_{1}\prec\ldots\prec
\gamma_{n}\prec\beta;\{m_{i}\}}C_{\{m_{i}\},\{\gamma_{i}\}}e_{\gamma_{1}}^{m_{1}}
e_{\gamma_{2}}^{m_2}\cdots e_{\gamma_{n}}^{m_n},
\end{eqnarray}
where $\sum_{i}^{n}m_{i}\gamma_{i}=\alpha+\beta$, and the coefficients $C_{\cdots}$ are
rational functions of $q$ and they do not depend on the Cartan elements $k_{\alpha_{i}}
,\: i=1,2,\ldots k$, and also
\begin{eqnarray}\label{QA10'}
[e_{\beta},e_{-\alpha}]\;=\;\sum
C_{\{m_{i}^{}\},\{\gamma_{i}^{}\};\{m'_{j}\},\{\gamma_{j}'\}}'
e_{-\gamma_1^{}}^{m_1}e_{-\gamma_{2}^{}}^{m_2}\cdots e_{-\gamma_p^{}}^{m_p}
e_{\gamma_1'}^{m_1'} e_{\gamma_2'}^{m_2'}\cdots e_{\gamma_{s}'}^{m_s'}
\end{eqnarray}
where the sum is taken on $\gamma_{1},\ldots,\gamma_{p},\gamma_{1}',\ldots,\gamma_{s}'$
and $m_{1}^{},\ldots,m_{p}^{}, m_{1}',\ldots,m_{s}'$ such that
\[
\gamma_{1}\prec\ldots\prec\gamma_{p}\prec\alpha\prec\beta\prec\gamma_{1}'\prec\ldots
\prec\gamma_{s}',\quad\sum_{l}(m_{l}'\gamma_{l}'-m_{l}\gamma_{l})=\beta-\alpha
\]
and the coefficients $C_{\cdots}'$ are rational functions of $q$ and $k_{\alpha}$ or
$k_{\beta}$. The monomials $e_{\gamma_{1}}^{n_{1}} e_{\gamma_{2}}^{n_{2}} \cdots
e_{\gamma_{p}}^{n_{p}}$ and $e_{-\gamma_{1}}^{n_{1}}e_{-\gamma_{2}}^{n_{2}}\cdots
e_{-\gamma_{p}}^{n_{p}}$, ($\gamma_{1} \prec\gamma_{2}\prec\cdots\prec\gamma_{p}$),
generate (as a linear space over $U_q(\mathcal{H})$) subalgebras
$U_{q}(\mathfrak{b}_{+})$ and $U_{q}(\mathfrak{b}_{-})$ correspondingly. The monomials
\begin{eqnarray}\label{QA10''}
e_{-\gamma_{1}^{}}^{n_{1}^{}}e_{-\gamma_{2}^{}}^{n_{2}^{}}\cdots
e_{-\gamma_{p}^{}}^{n_{p}^{}} e_{{\gamma}_{1}'}^{n_{1}'}e_{{\gamma}_{2}'}^{n_{2}'}\cdots
e_{{\gamma}_{s}'}^{n_{s}'},
\end{eqnarray}
where $\gamma_{1}\prec\gamma_{2}\prec\cdots\prec\gamma_{p}$ and $\gamma_{1}'\prec
\gamma_{2}'\prec\cdots\prec\gamma_{s}'$), generate $U_{q}(\mathfrak{g})$ over
$U_{q}(\mathcal{H})$.
\end{proposition}
Here the algebra $U_q(\mathcal{H})$ is generated by the Cartan elements $k_{\alpha_i}^{}$
($i=1,2,\ldots,r$).

Any formal Taylor series of $TU_q(\mathfrak{g})$ is constructed from the monomials of the
form (\ref{QA10''}) with coefficients, which are rational functions of the Cartan
elements $k_{\alpha_i}$ ($i=1,2,\ldots,r$), and nonnegative integers $n_{1}^{},
n_{2}^{},\ldots,n_{p}^{}, n_{1}',n_{2}',\ldots,n_{s}'$ are subjected to the constraints
of the type (\ref{aff8}) for all monomial of the given series.

By definition, the extremal projector for $U_{q}(\mathfrak{g})$ is a nonzero element
$p:=p(U_q(\mathfrak{g}))$ of the Taylor extension $TU_q(\mathfrak{g})$, $p\in
TU_q(\mathfrak{g})$, satisfying the equations
\begin{eqnarray}\label{QA11}
e_{\alpha_{i}}p=pe_{-\alpha_{i}}=0\quad(\forall\;\alpha_{i}\in\Pi), \qquad p^{2}=p.
\end{eqnarray}
We fix some normal ordering $\vec{\underline{\Delta}}_+(\mathfrak{g})$ and let
$\{e_{\pm\gamma}\}$ ($\gamma\in\underline{\Delta}_+(\mathfrak{g})$) be the corresponding
Cartan-Weyl generators. The following statement holds for any quantized
finite-dimensional contragredient Lie (super)algebra $g$ \cite{T3,KT1,T5}\footnote{The
theorem is also valid for the quantized infinite-dimensional affine Kac-Moody
(super)\-algebras, but in this case the formulas (\ref{QA13}) and (\ref{QA14}) for the
imaginary roots $\gamma=m\delta$ ($m\in\mathbb{N}$) should be more detailed (see
\cite{KT1, KLT} as examples).}.
\begin{theorem}
The equations (\ref{QA11}) have a unique nonzero solution in the space of the Taylor
extension $TU_{q}(\mathfrak{g})$ and this solution  has the form
\begin{eqnarray}\label{QA12}
p\;=\;\prod_{\gamma\in\vec{\underline{\Delta}}_+(\mathfrak{g})}p_{\gamma}^{},
\end{eqnarray}
where the order in the product coincides with the chosen normal ordering of
$\underline{\Delta}_+(\mathfrak{g})$ and the elements $p_{\gamma}$ are defined by the
formulae
\begin{eqnarray}\label{QA13}
p_{\gamma}^{}\;=\;\sum_{m\geq0}\frac{(-1)^{m}}{(m)_{\bar{q}_{\gamma}^{}}!}
\varphi_{\gamma,m}^{}e_{-\gamma}^{m}e_{\gamma}^{m},
\end{eqnarray}
\begin{eqnarray}\label{QA14}
\varphi_{\gamma,m}^{}\!\!&=\!\!&\frac{(q-q^{-1})^{m} q^{-\frac{1}{4}m(m-3)
(\gamma,\gamma)}q^{-m\rho(\gamma)}}{(a(\gamma))^m\prod\limits_{l=1}^{m}
\Bigl(k_{\gamma}^{}
q^{\rho(\gamma)+\frac{l}{2}(\gamma,\gamma)}-(-1)^{(l-1)\vartheta(\gamma)}k_{\gamma}^{-1}
q^{-\rho(\gamma)- \frac{l}{2}(\gamma,\gamma)}\Bigr)}.
\end{eqnarray}
Here $\rho$ is the linear function such that
$\rho(\alpha_{i})=\frac{1}{2}(\alpha_{i},\alpha_{i})$ for all simple roots
$\alpha_i\in\Pi$; $a(\gamma)$ is a factor in the relation (\ref{QA9});
$\bar{q}_{\gamma}^{}:=(-1)^{\vartheta(\gamma)}q^{-(\gamma,\gamma)}$; the symbol $(m)_q!$
is given as follows:
\begin{eqnarray}\label{QA15}
(m)_q!\;=\;(m)_q(m-1)_q\cdots(1)_q,\quad (0)_q!\;=\;1, \quad
(m)_{q}\;:=\;\frac{q^{m}-1}{q-1}.
\end{eqnarray}
\end{theorem}
In the limit $q\to1$ we obtain the extremal projector for the (super)algebra
$\mathfrak{g}$: $\lim\limits_{q\to1}p(U_q(\mathfrak{g}))=p(\mathfrak{g})$.

A proof of the theorem actually reduces to the proof for the case of the quantized
(super)algebras of rank 2, and it is similar to the case of non-deformed
finite-dimensional simple Lie algebras \cite{AST2,T4}.

\setcounter{equation}{0}
\section{Bibliographic comments on applications of\\
extremal projectors}

For convenience of the reader, we present here the most important references related to
the applications of extremal projectors.

(i) An explicit description of irreducible representations of Lie (super)algebras
(construction of various bases, action of generators, and their properties). Results
(Gelfand--Tsetlin bases) are contained in the following sources: \cite{AST3} for
$\mathfrak{su}(n)$, V.N. Tolstoy (1975, unpublished) for so(n), \cite{SST1}  for $G_2$,
\cite{M} for $\mathfrak{sp}(2n))$, \cite{BerTol} for $\mathfrak{osp}(1|2)$, \cite{TIS}
for $\mathfrak{gl}(n|m)$, \cite{T3,T5} for $U_q(\mathfrak{su}(n))$, \cite{PalTol} for
$U_q(\mathfrak{su}(n|1))$ and \cite{ABGST} for $U_q(\mathfrak{su}(n,1))$.

(ii) Theory of the Clebsch--Gordan coefficients of simple Lie algebras. Results can be
found in \cite{PlSmTo1}--\cite{PlSmTo3} for $\mathfrak{su}(3)$,
\cite{STKh1}--\cite{STKh4} for $U_q(\mathfrak{su}(2))$, \cite{AST4, AST5} for
$U_q(\mathfrak{su}(3))$ and \cite{TolDra} for $U_q(\mathfrak{su}(n))$.

(iii) Description of reduction algebras (Mikelson's algebras). Results are given in
\cite{Zh1} for $A_n$, $B_n$, $C_n$, and $D_n$; \cite{T7} for $\mathfrak{su}(n|m)$ and
$\mathfrak{osp}(m|2n)$; \cite{ABGST} for $U_q(\mathfrak{gl}(n))$.

(iv) Description of Verma modules of Lie (super)algebras (singular vectors and their
properties). Results for the simple Lie algebras can be found in \cite{Zh2}.

(v) Construction of solutions to the Yang--Baxter equation by using the projection
operators. Results for $\mathfrak{u}(3)$ and $\mathfrak{u}(n)$ are given in
\cite{ST1,TarVar}.

(vi) Relation between extremal projection operators and integral projection operators.
Results for simple Lie groups and algebras can be found in \cite{LezSav}.

(vii) Relation between extremal projection operators and canonical elements. Results for
$q$-boson Kashiwara algebras were obtained in \cite{Nak}.

(viii) Generalization of extremal projection operators. Results for $\mathfrak{sl}(2)$
and  $U_q(\mathfrak{sl}(2))$  can be found in \cite{T7} and \cite{DT}, respectively.

(ix) Construction of indecomposable representations. Results for $U_q(\mathfrak{sl}(2))$
are given in \cite{DT}.

\bibliographystyle{amsalphax} 
\end{document}